\documentclass[12pt]{article}
\usepackage{amsmath}
\usepackage{theorem}

\newtheorem{lemma}{Lemma}[section]

\newtheorem{theorem}[lemma]{Theorem}
\newtheorem{corollary}[lemma]{Corollary}
\newtheorem{definition}[lemma]{Definition}
\newtheorem{remark}[lemma]{Remark}

\newcommand{\N}{{\bf N}}

\newcommand{\R}{{\bf R}}
\newcommand{\C}{{\bf C}}

\newcommand{\cA}{{\cal A}}
\newcommand{\cB}{{\cal B}}

\newcommand{\cE}{{\cal E}}

\newcommand{\cint}{\iint\limits_{\C}}
\newcommand{\rint}{\int\limits_{-\infty}^{\infty}}

\newcommand{\bra}[1]{\left( #1 \right)}

\newcommand{\norm}[1]{\Vert #1 \Vert}

\newcommand{\Proof}{\underbar{Proof}{\hskip 0.1in}\par}

\newcommand{\Schrodinger}{Schr\"odinger }

\newcommand{\beq}{\begin{equation}}

\newcommand{\kik}[1]{C^{\infty}_{c}\bra{#1}}

\newcommand{\sz}{\langle z \rangle}

\newcommand{\sx}{\langle x \rangle}

\newcommand{\algebra}{\cA}

\begin{document}
\title{\bf Functional Calculus for Semi-Bounded Operators}
\author{Narinder S Claire}
\date{}
\maketitle
\begin {abstract}
We build on the work by Davies, extending the Helffer-Sj\"ostrand Functional Calculus domain for semi-bounded operators
on Banach spaces given a priori controlled growth of the resolvents. We employ Seeley's Extension Theorem to extend 
smooth functions on the half line to the whole line and thus indirectly define functions of these operators. 
\end{abstract}

\section*{Introduction}
Helffer and Sj\"ostrand \cite{hs} introduced a new formula into the field of spectral theory
\begin{equation}\label{e:hs}
f\bra{H} := -\frac{1}{\pi}\cint \frac{\partial \tilde{f}}{\partial \overline{z}}
\bra{z-H}^{-1}dxdy
\end{equation}
for self-adjoint operators. 
It was later shown by Davies \cite{Da1} that the Helffer-Sj\"ostrand formula had 
greater implications in spectral theory than it was originally intended for. He
\underline{constructed} a functional calculus for a much wider class of operators on spaces other than
Hilbert  
under certain assumptions on the norms of the resolvents. He in addition showed that
in the special case of self-adjoint operators
this coincided with the classical functional calculus.\\
We show that we can widen the class of functions when we consider
operators with spectrums which are bounded below, since positive operators take special
treatment in the analysis of partial differential equations. The results enable explicit tractable definitions of 
certain functions of partial differential operators in particular the Heat Semigroup.\\ 
Our hypothesis in the pending analysis is\par
{\em  $H$ is a closed densely defined operator on a Banach space $\cB$, with spectrum $\sigma\bra{H} \subset \R$.
It has resolvent operators $\bra{z-H}^{-1}$ defined and bounded for all $z \in \C$ satisfying:
\begin{equation}\label{e:hypothesis}
\norm{\bra{z-H}^{-1}} \leq c|Imz|^{-1}\bra{\frac{\langle z \rangle}{|Imz|}}^{\alpha}
\end{equation}
for some $\alpha \geq 0$ and all $z \not\in \R$, where $\langle z \rangle :=\bra{1+|z|^2}^{\frac{1}{2}}$} 

\section{Preliminaries}
We introduce our main concepts.
\subsection{Algebra of Slow Decreasing functions}

\begin{definition}
 $\sz := \bra{1+|z|^2}^{\frac{1}{2}}$
for all complex $z$  
\end{definition}
\begin{definition}
For $\beta \in \R$ let
$ S^{\beta}$ to be the set of all complex-valued smooth functions defined on $\R$ such that for every $n$
there is a positive constant
$c_n$ where $$|\frac{d^nf}{dx^n}|\, \leq\, c_n \sx^{\beta-n}$$
\end{definition}
\begin{definition}
We define the Algebra $\cA$ as,
\begin{equation}
\cA := \bigcup\limits_{\beta<0} S^{\beta}
\end{equation}
\end{definition}
\begin{remark}\label{rem:typ}
The strict inequality $\beta < 0$ is of importance and we conjecture that it
cannot be extended to
$\beta \leq 0$.
\end{remark}
\begin{lemma}[Davies \cite{Da2}]
$\cA$ is an algebra under pointwise multiplication.
If $f \in \cA$ then the expression
\begin{equation}
\norm{f}_{\cA_n} := \sum\limits_{r=0}^{n}\rint |\frac{d^{r}f}{dx^r}|\sx^{r-1} dx
\end{equation}
defines a norm on $\cA$ for each $n$
moreover $\kik{\R}$ is dense in $\cA$ with this norm.
The completion $\cA_n$ is a Banach space.
\end{lemma}
\begin{lemma}
The function $\sx^{\beta}$ is in $\cA$ for each $\beta <0$
\end{lemma}
\Proof
The statement follows from the observation that if $\beta <0 $ and $m\geq n$ then
$$ x^n\sx^{\beta-m}\leq \sx^{\beta}$$
and $$\frac{d\bra{x^n\sx^{\beta-m}}}{dx}=nx^{n-1}\sx^{\beta-m} +2\bra{\beta-m}x^{n+1}\sx^{\beta-m-1}$$

\begin{lemma}\label{lem:stable}
If $f\in\cA$ and $\phi \in S^0$ then $\phi f \in \cA$
\end{lemma}
\Proof
The statement follows from the inequality
\begin{eqnarray*}
|\frac{d^r\bra{\phi\bra{x}f\bra{x}}}{dx^r}| & =&|\sum\limits_{m=0}^rc_m\frac{d^{r-m}\bra{\phi\bra{x}}}{dx^{r-m}}\,
\frac{d^m\bra{f\bra{x}}}{dx^m}| \\
&\leq& c_r\sum\limits_{m=0}^r |\frac{d^{r-m}\bra{\phi}}{dx^{r-m}}|\,
|\frac{d\bra{f\bra{x}}^m}{d^mx}| \\
&\leq& c_{r,\phi}\sum\limits_{m=0}^r\,\sx^{\beta-r}\\
&\leq& c_{r,\phi}\sx^{\beta-r}
\end{eqnarray*}

\subsection{The Helffer Sj\"ostrand formula}
We introduce the concept of almost
analytic extensions due to H\"ormander \cite{ho}.
\begin{definition}
Let $\psi\bra{s}$ be a smooth function of compact support on $\R$
such that
\begin{equation*}
\psi\bra{s} : = \begin{cases} 1 & if |s|\leq 1 \\ 0 & if |s| \geq 2 \end{cases}
\end{equation*}
then we define 
\begin{equation}\label{e:please}
\sigma\bra{x,y} := \psi\bra{\frac{y}{\sx}}
\end{equation}
\end{definition}

\begin{definition}
Given $f \in \algebra$ we define an almost
analytic extension $\tilde{f}$ to the complex plane
\begin{equation}\label{e:ae}
\tilde{f}\bra{x,y} := \bra{\sum\limits_{r=0}^n \frac{d^{r}f\bra{x}}{dx^r}\frac{\bra{iy}^r}{r!}}\sigma\bra{x,y}
\end{equation}

and define
\begin{equation}\label{e:get}
\frac{\partial \tilde{f}}{\partial \overline{z}} := \frac{1}{2}\bra{\frac{\partial \tilde{f}}{\partial x}+
i\frac{\partial \tilde{f}}{\partial y}}
\end{equation}
\end{definition}
\begin{definition}
Given $f \in \algebra$ and $H$ satisfying our initial hypothesis
we define 
\begin{equation}\label{e:hs}
f\bra{H} := -\frac{1}{\pi}\cint \frac{\partial \tilde{f}}{\partial \overline{z}}
\bra{z-H}^{-1}dxdy
\end{equation}
\end{definition}
\par
We recall some important results from \cite{Da1} showing that we do indeed have a functional calculus
\begin{lemma}[Davies [\cite{Da1}]\label{e:important}
$ $  \begin{enumerate}
\item If $n>\alpha$ then subject to \eqref{e:hypothesis} the integral \eqref{e:hs} is norm
convergent for all $f$ in $\cA$ and
$$ \norm{f\bra{H}}\leq c\norm{f}_{n+1}$$
\item The operator $f\bra{H}$ is independent of $n$ and the cut-off function $\sigma$, subject to
$n>\alpha$\\
\item If $f$ is a smooth function of compact support disjoint from the spectrum of $H$ then $f\bra{H}=0$\\
\item If $f$ and $g$ are in $\cA$ then
$$\bra{fg}\bra{H}=f\bra{H}g\bra{H}$$
\item If $z \not\in \R$ and $g_z\bra{x} := \bra{z-x}^{-1}$ for real $x$ then
$g_z \in \cA$ and
$$g_z\bra{H}=\bra{z-H}^{-1}$$
\end{enumerate}
\end{lemma}
\section{Semi-bounded Operators}
We modify our main hypothesis by assuming the spectrum of $H$ is bounded below and without 
loss of generality $\sigma\bra{H} \subseteq [0,\infty)$.\\ 
We introduce a new ring
of functions $\cA^+$
\begin{definition} 
$S^{\beta}_+$ is the set of smooth functions on $\R^+\cup\{0\}$ with the same decaying property
as $S^\beta$ that is for every $n$ there is positive constant $c_n$ such that $$|\frac{d^nf}{dx^n}|\,
\leq\, c_n \sx^{\beta-n}$$
Then $\cA^+$ is defined appropriately and similarly we define the Banach space $\cA^+_n$ with norm
\begin{equation}
\norm{f}_{\cA^+_n} := \sum\limits_{r=0}^{n}\int\limits_0^{\infty} |\frac{d^{r}f}{dx^r}|\sx^{r-1} dx
\end{equation}
\end{definition}
\subsection{Seeley's Extension Theorem}
We present a theorem due to Seeley \cite{se} which gives a linear extension operator
for smooth functions from the half space to the whole space. This extension operator is continuous
for many topologies including uniform convergence of each derivative. We demonstrate a brief
proof as it contributes to the proof of continuity for our topology $\norm{ }_{\cA^+_n}$ for each n.
\begin{definition}
Given  $f \in \cA^+$, $\phi \in \cA$ and real $a$ we define two operators on $\cA^+$,
$$\bra{T_af}\bra{x}=f\bra{ax}$$
$$\bra{S_{\phi}f}\bra{x}=\phi\bra{x}f\bra{x}$$
\end{definition}
\begin{theorem}Seeley's Extension Theorem.\\
There is a linear extension operator
$$ \cE : C^{\infty}[0,\infty) \longrightarrow C^{\infty}\bra{\R}$$
such that for all $x>0$
$$\bra{\cE f}\bra{x}=f\bra{x} $$
\end{theorem}
The proof of the theorem centres on the following lemma
\begin{lemma}\label{lem:seeley}
There are sequences $\{a_k\},\, \{b_k\}$ such that
\begin{enumerate}
\item $b_k <0$
\item $\sum\limits_{k=0}^{\infty}|a_k||b_k|^n<\infty$\label{p:two}
 for all non-negative integers $n$ \label{p:two}
\item $\sum\limits_{k=0}^{\infty}a_k \bra{b_k}^n =1$ for all non-negative integers $n$
\item $b_k \rightarrow -\infty$
\end{enumerate}
\end{lemma}
\Proof
See \cite{se}. We recall from the proof $b_k=-\bra{2^k}$ 
 and $|a_k|<e^42^{-\frac{k^2-3k}{2}}$
\par
Proof of theorem \newline
Let $\phi \in \kik{\R}$ such that
$$\phi\bra{x}=\begin{cases} 1 & x \in [0,1] \\ 0 & x \geq 2 \\ 0 & x \leq -1 \end{cases} $$
we construct $\cE$ 
$$ \bra{\cE f}\bra{x} := \begin{cases} \sum\limits_{k=0}^\infty a_k \bra{T_{b_k}S_{\phi}f}\bra{x} &x <0\\
                                                f\bra{x} & x\geq 0
\end{cases} $$

The series is convergent since for all negative $x$ the sum has only finite non-zero terms.
It is evident that 
$ \bra{\cE f}\bra{0}= f\bra{0}$
and when $N\in\N $ then for all $x>-\frac{1}{2^{N}}$ $$\phi^{\bra{n}}\bra{x}=0$$ for all positive $n$
hence
$$\bra{\cE f}^{\bra{n}}\bra{x} =  
\sum\limits_{k=0}^{\infty} a_kb_k^n\phi\bra{b_k x}f^{\bra{n}}\bra{b_k x}$$
and we deduce that for all n $\lim\limits_{x\rightarrow 0-}\frac{d^n\bra{\cE f}}{dx^n}\bra{x}=\lim\limits_{x\rightarrow 0+}
\frac{d^n\bra{\cE f}}{dx^n}\bra{x}$ to complete the proof.
\begin{lemma}
If  $a>1$ then $\norm{T_a}_{\cA^+_n \rightarrow \cA^+_n}\leq a^n$
\end{lemma}
\Proof
\begin{eqnarray*}
\norm{T_af}_{\cA^+_n}&=&\sum\limits_{r=1}^n\int\limits_0^{\infty} |\frac{d^rf\bra{ax}}{dx^r}|\sx^{r-1}dx\\
 & \leq &\sum\limits_{r=1}^na^r\int\limits_0^{\infty} |\frac{d^rf\bra{ax}}{d\bra{ax}^r}|\langle ax
\rangle^{r-1}d\bra{ax}\\
& = & \sum\limits_{r=1}^na^r\int\limits_0^{\infty} |\frac{d^rf\bra{x}}{dx^r}|\sx^{r-1}dx
\end{eqnarray*}
and the inequality follows.
\begin{lemma}If $\phi\in\cA$ then 
$S_{\phi}$ is a bounded operator with respect
 to each norm $\norm{\,}_{\cA^+_n}$
\end{lemma}
\Proof
A simple application of Leibnitz gives
$$ \frac{d^r\bra{\phi\bra{x}f\bra{x}}}{dx^r}=\sum\limits_{m=0}^rc_r\frac{d^{r-m}\bra{\phi\bra{x}}}{dx^{r-m}}\,
\frac{d^m\bra{f\bra{x}}}{dx^m}$$ then
\begin{eqnarray*}
|\frac{d^r\bra{\phi\bra{x}f\bra{x}}}{dx^r}|&\leq &c_r\sum\limits_{m=0}^rd_{r-m,\phi}\,\sx^{\beta-\bra{r-m}}
\frac{d^m\bra{f\bra{x}}}{dx^m}\\
&\leq & c_{r,\phi}\sum\limits_{m=0}^r\,\sx^{m-r}\frac{d^m\bra{f\bra{x}}}{dx^m}
\end{eqnarray*}
and so we integrate to give
\begin{eqnarray*}
\int\limits_{0}^{\infty}
|\frac{d^r\bra{\phi\bra{x}f\bra{x}}}{dx^r}|\sx^{r-1}dx & \leq &
c_{r,\phi}\sum\limits_{m=0}^r\int\limits_{0}^{\infty}|\frac{d^m\bra{f\bra{x}}}{dx^m}|\sx^{m-1}dx\\
& = & c_{r,\phi}\norm{f}_{\cA^+_r}
\end{eqnarray*}
and hence we have our estimate
\begin{eqnarray*}
\norm{S_{\phi}f}_n & = & \sum\limits_{r=0}^n\int\limits_{0}^{\infty}
|\frac{d\bra{\phi\bra{x}f\bra{x}}^r}{d^rx}|\sx^{r-1}dx\\
& \leq & c_{n,\phi}\sum\limits_{r=0}^n \norm{f}_{\cA^+_r} \\
&\leq &c_{n,\phi}\norm{f}_{\cA^+_n}
\end{eqnarray*}
\begin{theorem}
Seeley's Extension Operator is a bounded operator on each of the normed vector spaces $\cA_n^+$
\end{theorem}
\Proof
\begin{eqnarray*}
\norm{\cE f}_{\cA_n} &=& \sum\limits_{r=0}^n \int\limits_{-\infty}^{\infty}|\frac{d^r\bra{\cE f}}{dx^r}|\sx^{r-1}dx\\
&=&\sum\limits_{r=0}^n \int\limits_0^{\infty}|\frac{d^rf\bra{x}}{dx^r}|\sx^{r-1}dx +
\sum\limits_{r=0}^n \int\limits_{-\infty}^{0}|\sum\limits_0^{\infty}a_k\frac{d^r\bra{\phi\bra{b_kx}f\bra{b_kx}}}{dx^r}|\sx^{r-1}dx
\\ &=&\norm{f}_{\cA^+_n} + \norm{\sum\limits_{k=0}^{\infty}a_k T_{-b_k}S_{\phi}f}_{\cA^+_n}\\
& \leq & \norm{f}_{\cA^+_n} +\sum\limits_{k=0}^{\infty}|a_k|\,\norm{S_{\phi}}\,\norm{|T_{-b_k}}\norm{f}_{\cA^+_n}\\
& \leq &\norm{f}_{\cA^+_n} +\bra{\sum\limits_{k=0}^{\infty}|a_k|\,|b_k|^n}\,c_{n,\phi}\norm{f}_{\cA^+_n}
\end{eqnarray*}
and hence the extension operator is continuous.

\subsection{The Functional Calculus}
\begin{remark}
If $f$ and $g$ are elements of $\cA$ such that $f|_{[0,\infty]}=g|_{[0,\infty]}$ and the spectrum of $H$ is $[0,\infty)$
then it is not necessary that $supp\bra{f-g}\cap\sigma\bra{H}$ is empty \\since $supp\bra{f-g}\cap\sigma\bra{H}=\{0\}$
is possible and lemma \ref{e:important} cannot be applied. This renders our problem non-trivial and 
justifies the technical detour. 
\end{remark}
\begin{theorem}
If $f$ is a smooth function on $\R$ of compact support such that
$$supp\bra{f}=[-a,0]$$ and $H$ is an operator satisfying our modified
hypothesis with $\sigma\bra{H}\subseteq [0,\infty]$
then $$f\bra{H}=0$$
\end{theorem}
\Proof
Let $\epsilon \in \bra{0,1}$ and define $$f_{\epsilon}\bra{x}:=f\bra{x+\epsilon}$$
so that $supp\bra{f_{\epsilon}}=[-\bra{a+\epsilon},-\epsilon]$.\\
By Lemma \ref{e:important} $ \bra{f_{\epsilon}\bra{H}}=0$.
For all $n$ there are constants $c_n\geq 0$ such that 
$$\norm{\frac{d^nf}{dx^n}-\frac{d^nf_{\epsilon}}{dx^n}}_{\infty}\leq c_n\epsilon$$
then 
\begin{eqnarray*}
\norm{f\bra{H}}&=&\norm{f\bra{H}-f_{\epsilon}\bra{H}}\\
&\leq & \sum\limits_{r=0}^n \int\limits_{-\bra{a+1}}^0
|\frac{d^rf\bra{x}}{dx^r}-\frac{d^rf_{\epsilon}\bra{x}}{dx^r}|\sx^{r-1}dx\\
&\leq &\sum\limits_{r=0}^n \epsilon c_r \int\limits_{-\bra{a+1}}^0 \sx^{r-1}dx\\
&=& \epsilon k_{n,f}
\end{eqnarray*}
hence our result.
\begin{corollary}\label{last}
If $f$ and $g$ are in $\cA$ such that $f|_{[0,\infty]}=g|_{[0,\infty]}$
and $\sigma\bra{H}\subseteq [0,\infty]$ then $f\bra{H}-g\bra{H}=0$
\end{corollary}
\begin{theorem}
If $H$ is a closed densely defined operator on a Banach space $\cB$, with spectrum $\sigma\bra{H} \subset [0,\infty)$.
with resolvent operators $\bra{z-H}^{-1}$ defined and bounded for all $z \in \C$ satisfying:
\begin{equation}
\norm{\bra{z-H}^{-1}} \leq c|Imz|^{-1}\bra{\frac{\langle z \rangle}{|Imz|}}^{\alpha}
\end{equation}
for some $\alpha \geq 0$ and all $z \not\in \R$\\
then there is a functional calculus $\gamma_H :\cA^+\rightarrow {\bf B}\bra{\cB}$ 
such that for all $f\in \cA^+\bigcap\cA$ $$\gamma_H\bra{f}=-\frac{1}{\pi}\cint \frac{\partial \tilde{f}}{\partial 
\overline{z}}
\bra{z-H}^{-1}dxdy$$
\end{theorem}
\Proof
Let $f^+ \in \cA^+$, then by Seeley's Extension Theorem there exists an extension $f\in\cA$.
We define $\gamma_H\bra{f^+}:=f\bra{H}$. This definition is independent of the particular 
extension by corollary \ref{last}.The functional analytic properties are inherited from the extension.
   
\section*{Acknowledgements}
This research was funded by an EPSRC Ph.D grant 95-98 and my mother. I would like to thank Brian Davies
for giving me this problem and his encouragement since. I am grateful to all the research students
in the Mathematics Dept. Kings College, London 1995-2000 for all their help and support. I thank Mark Owen 
and Colin Mason for some helpful discussions I am indebted 
 to Anita for all her support.

\vskip 0.3in
Department of Mathematics \newline
Strand \newline
London WC2R 2LS \newline
King's College \newline
England \\
e-mail: nclaire@mth.kcl.ac.uk
\vfill

\begin{thebibliography}{99}
\bibitem{Da1} E.B. Davies, {\em The Functional Calculus},
\newline J. London
Math. Soc (2) 52 (1995) 166-176
\bibitem{Da2} E.B. Davies, {\em Spectral Theory and Differential
Operators},\newline C.U.P 1995
\bibitem{hs} B. Helffer, J. Sj\"ostrand. {\em Equation de \Schrodinger avec Champ Magnetique et
Equation de Harper}, 
\newline Lexture Notes in Physics 345 \Schrodinger Operators  Springer (1989) p118-197
\bibitem{ho} L. H\"ormander, {\em Linear Partial Differential
Operators},\newline Springer 1993  
\bibitem{se} S.T. Seeley {\em Extensions of $C^\infty$ functions defined on a half
space}, \newline Proc. Amer. Math. Soc 15 1964
\end{thebibliography}
\end{document}